\documentstyle[11pt,amscd]{amsart}
\newcommand{\proj}{{\mathbf{P}}}
\newcommand{\com}{{\mathbb{C}}}
\newcommand{\rarr}{\rightarrow}

\begin{document}
\title{Convex rationally connected varieties}
\author{R. Pandharipande}
\date{7 January 2004}
\maketitle
\pagestyle{plain}
\baselineskip=16pt 

\setcounter{section}{-1}

\section{Introduction}
Let $X$ be a nonsingular projective variety over $\mathbb{C}$.
A morphism,
$$\mu: \proj^1 \rarr X,$$
is {\em unobstructed} if $H^1(\proj^1, \mu^* T_X)=0$.
The variety $X$ is {\em convex} if all morphisms $\mu: \proj^1 \rarr X$ are 
unobstructed.

A {\em rational curve} in $X$ is the image of 
a morphism 
$$\mu: \proj^1 \rarr X.$$
The variety $X$ is {\em rationally connected} if all pairs of points of $X$
are connected by rational curves.

Homogeneous spaces ${\mathbf G}/{\mathbf P}$ for connected
linear algebraic groups are
convex, rationally connected, nonsingular, projective varieties. 
Convexity is a consequence of the
global generation of the tangent bundle of ${\mathbf G}/{\mathbf P}$.
Rational connectedness is consequence of the 
rationality of ${\mathbf G}/ {\mathbf P}$.

The following speculation arose at dinner after an algebraic
geometry seminar at Princeton in the fall of 2002.

\noindent {\bf Speculation.} If $X$ is convex and rationally connected,
then $X$ is a homogeneous space.

\noindent The failure of the
speculation would perhaps be more interesting than the success.

\section{Complete intersections}

The only real evidence known to the author is the following
result for complete intersections in projective space.

\noindent {\bf Theorem.} {\em If $X \subset \proj^n$ is a
convex, rationally connected, nonsingular complete intersection, then
$X$ is a homogeneous space.}

\noindent{\em Proof.}
We first consider nonsingular complete intersections
of dimension 
at most 1:
\begin{enumerate}
\item[(i)]
in dimension 0, only points are rationally connected,
\item[(ii)]
in dimension 1, only  
$\proj^1$ is rationally connected.
\end{enumerate}
Hence, the rationally connected complete intersections of
dimension at most 1 are simply connected.

Let $X\subset \proj^n$ be a {\em generic} complete intersection of 
type $(d_1, \ldots, d_l)$.
Let $d= \sum_{i=1}^l d_i$. By the results of \cite{fano},
$X$ is rationally connected if and only if
$d \leq n$. Moreover, if $X$ is rationally connected, then
$X$ must be simply connected: if the dimension of 
$X$ at least 2, $X$ is simply connected by Lefschetz, see \cite{mil}.
 
Let $M$ denote the parameter space of lines in $X$. 
$M$ is a non-empty, nonsingular variety of dimension
$2n-2-d-l$. Non-emptiness can be seen by several methods. For example,
the nonvanishing in degree 1 of the
1-point series of the quantum cohomology of $X$ implies $M$ is
non-empty, see \cite{giv}, \cite{bour}. Nonsingularity is
a consequence of the genericity of $X$.
Let $\pi: U \rarr M$ denote the universal family of
lines over $M$, and let
$$\nu: U \rarr X$$
denote the universal morphism.

Let $L$ be a line on $X$. If $X$ is convex, the normal bundle $N_L$
of $L$ in $X$ must be semi-positive. If $N_L$ has a negative
line summand, then every
double cover of $L$, $$\mu: \proj^1 \rarr L \subset X,$$
is obstructed.
Since the degree of $N_L$ is $n-d-1$, we may assume
$d \leq n-1$.

Every semi-positive bundle on $\proj^1$ is generated by global
sections. Hence, 
if every line $L$ has semi-positive normal bundle, we easily conclude
the morphism $\nu$ is smooth and surjects onto $X$. 
The fiber of $\nu$ over $x\in X$
is the parameter space of lines passing through $x$.

We now consider the Leray spectral sequence for the
fibration $\nu$, see \cite{gh}. The Leray spectral sequence degenerates at
the $E_2$ term,
$$E_2^{pq} = H^p(X, R^q {\nu_*} \com).$$
Since $X$ is simply connected, all local systems
on $X$ are constant. Hence,
$$E_2^{pq} = H^p(X, R^q {\nu_*} \com)= H^p(X,\com) \otimes H^q(F,\com),$$
where $F$ denotes the fiber of $\nu$.

Let $p_U(t)$, $p_F(t)$, and $p_X(t)$ denote the Poincar\'e polynomials
of the manifolds $U$, $F$, and $X$.
We conclude,
$$p_U= p_F \cdot p_X.$$
On the other hand, since $U$ is a locally trivial fibration over
$M$, the polynomial
$p_{\proj^1}$ must divide  $p_U$.
Since $$p_{\proj^1}=1+t^2$$ is irreducible over the integers, 
we find $1+t^2$ divides either $p_F$ or $p_X$.

We have proven the following result. Let $X\subset \proj^n$ 
be a generic
complete intersection of type $(d_1, \ldots,d_l)$
satisfying $d\leq n-1$.
If every line of $X$ has a semi-positive normal bundle, then
either $p_F(i)=0$ or
$p_X(i)=0$.

Consider the fiber $F$ of $\nu$ over $x$. 
The dimension of $F$ is $n-1-d$. In fact,
$F$ is a complete intersection of type
$$(1,2,3,\ldots,d_1, 1,2,3,\ldots,d_2,\ldots, 1,2,3, \ldots, d_l)$$
in the projective space 
$\proj^{n-1}$ of lines of $\proj^n$ passing through
$x$.

If $p_F(i)=0$, then the type of $F$ must be one of the
three types allowed by the Lemma below.
If $p_X(i)=0$, then the type of $X$ must be one of the
three allowed by the Lemma.
Since, one of the two polynomial evaluations must vanish,
we conclude the type of $X$ must be either
$(1, \ldots,1)$  or $(1,\ldots,1,2).$
Clearly both are types of homogeneous varieties.

If $X$ is not of one of the two above types, then $X$ must contain a 
line $L$ for which $N_L$ has a negative line summand.
Since $X$ was assumed to be general, every nonsingular
complete intersection $Y$ of the type of $X$ must also contain
such a line by taking a limit of $L$.

Therefore, if the type of a nonsingular complete intersection $Y$ is
not $(1, \ldots,1)$  or $(1,\ldots,1,2)$, 
then $Y$ is not a convex, rationally connected
variety.
\qed

The proof of the Theorem also shows homogeneous complete
intersections in projective space must be of 
type
$(1, \ldots,1)$  or $(1,\ldots,1,2)$.

\noindent{\bf Lemma.} {\em  Let $Y \subset \proj^n$ be a nonsingular
complete intersection of dimension $k$.
Let $p_Y(t)$ be the Poincar\'e polynomial of $Y$.
If $p_Y(i)=0$, then 
one of the following three possibilities hold:
\begin{enumerate}
\item[(i)] the type of $Y$ is $(1,\ldots,1)$ and $k$ is
odd,
\item[(ii)] the type of $Y$ is
$(1,\ldots,1,2)$ and $k$ is odd,
\item[(iii)] the type of $Y$ is
$(1,\ldots,1,2)$, and $k= 2 \mod 4$.
\end{enumerate}}

\noindent{\em Proof.}
Let $Y\subset \proj^n$ be 
a nonsingular complete intersection of dimension $k$.
The cohomology of $Y$ is determined by the Lefschetz isomorphism
except in the middle (real) dimension $k$.
The cohomology determined by Lefschetz is of rank 1 in all
even (real) dimensions.
If $k$ is odd, then 
$$p_Y(t)= \sum_{q=0}^k t^{2q} + b_k t^k,$$
where $b_k$ is the $k^{th}$ Betti number.
We see $p_Y(i)=0$ if and only if $b_k=0$.
If $k$ is even, 
$$p_Y(t)= \sum_{q=0}^k t^{2q} + (b_k-1) t^k.$$
We see $p_Y(i)=0$ if and only if 
$k=2 \mod 4$   
and $b_k-1=1$.

Assume $p_Y(i)=0$.
Let $(e_1, \ldots, e_k)$ be the type of $Y$. Let $e$ be the
largest element of the type.

Let $Z \subset \proj^n$ be a nonsingular projective variety of dimension
$r$. Let $H_d\subset \proj^n$ be a general hypersurface of degree $d$.
The dimension,   
$$h^{r-1}(Z \cap H_d, \com),$$
is a non-decreasing function of $d$, see \cite{katz}.
Hence, we can bound $b_k$ for $Y$ from below by the middle
cohomology $b_k'$ of the complete intersection $Y'\subset \proj^n$ of type
$(e,1,\ldots,1)$,
$$b_k \geq b_k'.$$
 The variety $Y'$ may then be viewed
as a hypersurface of degree $e$ in the smaller projective space
$\proj^{k+1}$.

For a hypersurface $Y'\subset \proj^{k+1}$ of degree $e$, the
middle cohomology $b'_k$ is given by the following formula:
$$b_k'- \delta_{k}= \frac{(e-1)}{e}((e-1)^{k+1} - (-1)^{k+1}),$$
where $\delta_k$ is 1 if $k$ is even and 0 if $k$ is odd.
If $k$ is odd,
$$b_k' =  \frac{(e-1)}{e}((e-1)^{k+1} - 1).$$
Then, $b_k' >0$ if $e>2$.
If $k$ is even, 
$$b_k'-1 =  \frac{(e-1)}{e}((e-1)^{k+1} + 1).$$
Then, $b_k-1>1$ if $e>2$.
Therefore, we conclude $e \leq 2$.

If $e=1$, then case (i) of the Lemma is obtained. It is easy to
check $k$ must be odd for $p_Y(i)=0$ to hold.

Let $e=2$. If $Y$ is of type $(1,\ldots,1,2)$, then either case (ii)
or (iii)
of the Lemma is obtained. If $k$ is even,
$$k= 2\ \text{mod}\ 4$$ must be satisfied in order
for $p_Y(i)=0$ to hold.

If $Y$ is not of type $(1,\ldots,1,2)$, then
the next largest type of $Y$ must be at least 2.
As before, we may bound $b_k$ from below by the middle cohomology
$b_k'$ of the complete intersection of type $(2,2)$ in
$\proj^{k+2}$.
If $k$ is odd, the calculation below shows
$$b_k'= k+1 >0.$$
If $k$ is even, the calculation below shows
$$b_k'-1 = k+3>1.$$
In fact, the type of $Y$ {\em can not} contain two elements
greater than 1 in $p_Y(i)=0$.

Let $k\geq 0$. The Euler characteristic $\chi_{22}(k)$ 
of a nonsingular
complete intersection of type $(2,2)$ in $\proj^{k+2}$ is:
\begin{eqnarray*}
\int_{\proj^{k+2}} \big( \frac{2H}{1+2H} \big)^2 (1+H)^{k+3}
& = & \sum_{i=0}^k 2^{k+2-i}(-1)^{k-i}(k+1-i)\binom{k+3}{i}. \\
\end{eqnarray*}  
On the other hand, since
$$\hspace{-100 pt}
(k+3)(t-1)^{k+2} - (t-1)^{k+3} + (-1)^{k+3} = 
$$
$$\hspace{+100 pt}
\sum_{i=0}^{k+2} t^{k+2-i}(-1)^i (k+3-i-t)\binom{k+3}{i}, $$
we find:
\begin{eqnarray*}
(-1)^k\chi_{22}(k) =
k+2 +(-1)^{k}(k+2).
\end{eqnarray*}
The formulas for $b_k'$ then follow easily.
\qed

\section{Homogeneous complete intersections}
It is interesting to see how the homogeneous complete intersections
survive the above analysis.

First, consider a complete intersection $X\subset \proj^n$
of type $(1,\ldots,1)$.
Then, $F$ is of dimension $n-1-l$, and $X$ is a dimension $n-l$.
Both are complete intersections of hyperplanes.
Since one of $n-1-l$ and $n-l$ is odd, exactly one of
the conditions $p_F(i)=0$ or
$p_X(i)=0$ holds by part (i) of the Lemma.

Next, consider a complete intersection $X\subset \proj^n$
of type $(1,\ldots,1,2)$ where $l+1 \leq n-1$
Then, $F$ is of dimension $n-1-l-1$, and $X$ is of dimension
$n-l$. There are two cases:
\begin{enumerate}
\item[(i)] If $n-l-2$ and $n-l$ are odd, then both $p_F(i)=0$
and $p_X(i)=0$ by part (ii) of the Lemma.
\item[(ii)] If $n-l-2$ and $n-l$ are even, then one of
$(n-l-2)/2$ and $(n-l)/2$ is odd. Hence, exactly  one of the conditions 
$p_F(i)=0$ and $p_X(i)=0$ holds 
by part (iii) of the Lemma.
\end{enumerate}

\noindent Department of Mathematics\\
\noindent Princeton University\\
\noindent Princeton, NJ 08540

\end{document}